\documentclass[Svgc,gcfinal]{Svgc}
\usepackage{graphics}
\usepackage{amsmath}
\usepackage{amssymb}
\usepackage{enumerate}
\journalname{Graphs and Combinatorics}

\newcommand{\A}{{\cal A}}
\newcommand{\B}{{\mathcal B}}
\newcommand{\D}{{\sf D}}

\newcommand{\G}{{\cal G}}
\newcommand{\I}{{\mathcal I}}

\renewcommand{\S}{{\cal S}}

\newcommand{\vecx}{{\bf x}}

\newcommand{\Int}{{\rm Int}}
\newcommand{\supp}{{\rm supp}}

\begin{document}
\title{The Fine Intersection Problem for Steiner Triple Systems}
\subtitle{}
\author{Yeow Meng Chee\inst{1,2} \and Alan C. H. Ling\inst{3} \and Hao Shen\inst{4}
}                     


%
%
\institute{Division of Mathematical Sciences,
School of Physical and Mathematical Sciences,
Nanyang Technological University,
Singapore 637616, \email{ymchee@ntu.edu.sg}. \and
Card View Pte. Ltd.,
41 Science Park Road,
\#04-08A The Gemini,
Singapore Science Park II,
Singapore 117610, \email{ymchee@alumni.uwaterloo.ca}\and 
Department of Computer Science,
University of Vermont,
Burlington, Vermont 05405.
USA, \email{aling@emba.uvm.edu}. \and
Department of Mathematics,
Shanghai Jiao Tong University,
Shanghai 200030,
People's Republic of China, \email{haoshen@sjtu.edu.cn}.
}
\maketitle
\begin{abstract}
The intersection of two Steiner triple systems $(X,\A)$ and $(X,\B)$ is the set $\A\cap\B$.
The fine intersection problem for Steiner triple systems is to determine for each $v$, the
set $I(v)$, consisting of all possible pairs $(m,n)$ such that there exist two Steiner triple
systems of order $v$ whose intersection $\I$ satisfies $|\cup_{A\in\I} A|=m$ and
$|\I|=n$. We show that for $v\equiv 1$ or $3\pmod{6}$, $|I(v)|=\Theta(v^3)$, where
previous results only imply that $|I(v)|=\Omega(v^2)$.
\end{abstract}
%


\section{Introduction}
For a set $X$ and non-negative integer $k$, denote by $X\choose k$ the set of all
$k$-subsets of $X$. The {\em support} of $\A\subseteq 2^X$, denoted by $\supp(A)$,
is the set $\cup_{A\in\A} A$. A {\em set system} is a pair $(X,\A)$, where $X$ is a finite
set of {\em points}, and $\A\subseteq 2^X$. The elements of $A$ are called {\em blocks}.
The {\em order} of a set system is the number of points in the set system. Let $K$ be a set
of positive integers. The set $K$ is a {\em set of block sizes} for $(X,\A)$ if $|A|\in K$ for all
$A\in\A$.
A set system $(X,\A)$ is said to be $k$-{\em uniform} if $\A\subseteq{X\choose k}$.

Let $(X,\A)$ be a set system and let $\G=\{G_1,\ldots,G_s\}$ be a partition of $X$ into subsets
called {\em groups}. The triple $(X,\G,\A)$ is a {\em group divisible design} (GDD) when every
2-subset of $X$ is either contained in exactly one block or in exactly one group. We denote a
GDD $(X,\G,\A)$ by $K$-GDD if $K$ is a set of block sizes for $(X,\A)$. The {\em type} of a GDD
$(X,\G,\A)$ is the multiset $[|G| : G\in\G]$. When more convenient, we use the exponentiation
notation to describe the type of a GDD: a GDD of type $g_1^{t_1}\cdots g_s^{t_s}$ is a 
GDD where there are exactly $t_i$ groups of size $g_i$, $1\leq i\leq s$.

A $\{3\}$-GDD of type $1^v$ is a {\em Steiner triple system} (STS) of order $v$, and is denoted
STS$(v)$. It is well-known that an STS$(v)$ exists if and only if $v\equiv 1$ or $3\pmod{6}$
(see, for example, \cite{ColbournRosa:1999}). A {\em partial triple system} is a 3-uniform set
system $(X,\A)$ where every 2-subset of $X$ is contained in at most one block of $\A$.

The {\em intersection} of two $K$-GDDs (of the same type) $\D_1=(X,\G,\A_1)$
and $\D_2=(X,\G,\A_2)$ is the set $\I(\D_1,\D_2)=\A_1\cap\A_2$. $\D_1$ and $\D_2$ are said
to be {\em disjoint} if $\I(\D_1,\D_2)=\varnothing$. Let $f:2^{X\choose 3}\rightarrow \Gamma$.
The $f$-{\em intersection problem for $K$-GDDs of type $T$} is to determine the set
\begin{align*}
   & \Int_f(K,T) \\
= & \{r\in\Gamma : \text{$\exists$ two $K$-GDDs of type $T$, $\D_1$ and $\D_2$, with
$f(\I(\D_1,\D_2))=r$}\}.
\end{align*}
The interest in this paper is the $f$-intersection problem for STS in the case when $f=\Phi$,
where $\Phi:2^{X\choose 3}\rightarrow \mathbb{Z}^2_{\geq 0}$ is defined as follows:
\begin{eqnarray*}
\Phi(\S) = (|\supp(\S)|,|\S|), & & \text{for $\S\subseteq{X\choose 3}$.}
\end{eqnarray*}
We call this the {\em fine intersection problem} for STS for the reason that both the number of
blocks and the number of underlying points in the intersection are to be determined simultaneously.

All previous work on the intersection of STS can be cast in the context of $f$-intersection problems
for STS, for appropriate choices of $f$.

\begin{example}[Lindner and Rosa \cite{LindnerRosa:1975a}]
The classical problem of determining the possible number of blocks in the intersection of two STS$(v)$
is equivalent to the $f$-intersection problem for STS with $f(\S)=|\S|$.
\end{example}

\begin{example}[Hoffman and Lindner \cite{HoffmanLindner:1987}]
The flower intersection problem for STS$(v)$ is equivalent to the $f$-intersection problem for STS, with
\begin{eqnarray*}
f(\S)=\begin{cases}
|\S|,&\text{if $\S$ contains a set of $\frac{v-1}{2}$ blocks intersecting in a common point;} \\
\infty,&\text{otherwise.}
\end{cases}
\end{eqnarray*}
\end{example}

\begin{example}[Chee \cite{Chee:2004}]
The disjoint intersection problem for STS is equivalent to the $f$-intersection problem for STS, with
\begin{eqnarray*}
f(\S)=\begin{cases}
|\S|,&\text{if the blocks in $\S$ are pairwise disjoint;}\\
\infty,&\text{otherwise.}
\end{cases}
\end{eqnarray*}
This is also equivalent to determining the possible values of $n$ for which
$(3n,n)\in\Int_\Phi(\{3\},1^v)$.
\end{example}

The purpose of this paper is to initiate the study on the fine intersection problem for STS.

\section{Admissible Region}

In this section, we determine an admissible region for the fine intersection problem for STS; that is,
we determine a subset of $\mathbb{Z}^2_{\geq 0}$ which contains $\Int_\Phi(\{3\},1^v)$.

A partial triple system $(X,\A)$ of order $v$ is said to be {\em maximum} if for every partial
triple system $(X,\B)$ of order $v$, we have $|\B|\leq |\A|$. The number of blocks in a maximum
partial triple system is denoted $D(2,3,v)$. For a vector $\vecx=(x_1,\ldots,x_d)\in\mathbb{Z}^d$,
we denote by $\vecx|_i$, $1\leq i\leq d$, the value $x_i$. For a set of vectors
$S\subseteq\mathbb{Z}^d$, we denote by $S|_i$, $1\leq i\leq d$, the {\em projection} of $S$
on the $i^{\rm th}$ dimension: $S|_i=\{\vecx|_i : \vecx \in S\}$.

We state below some prior results that are useful in establishing an admissible region for
$\Int_\Phi(\{3\},1^v)$.

\begin{theorem}[Lindner and Rosa \cite{LindnerRosa:1975a}]
\label{LR}
Let $b(v)=v(v-1)/6$. Then for all $v\equiv 1$ or $3\pmod{6}$, $v\not=9$,
$\Int_{|\cdot|}(\{3\},1^v)=\{0,\ldots,b(v)\}\setminus\{b(v)-5,b(v)-3,b(v)-2,b(v)-1\}$, and
$\Int_{|\cdot|}(\{3\},1^9)=\{0,1,2,3,4,6,12\}$.
\end{theorem}

\begin{theorem}[Doyen and Wilson \cite{DoyenWilson:1973}]
\label{DW}
Let $v,w\equiv 1$ or $3\pmod{6}$, and $w<v$. There exists an STS$(v)$ containing an
STS$(w)$ if and only if $v\geq 2w+1$.
\end{theorem}

\begin{lemma}
\label{admissible}
Let $(m,n)\in\Int_\Phi(\{3\},1^v)$. Then the following conditions hold:
\begin{enumerate}[{\rm (i)}]
\item $m/3\leq n\leq D(2,3,m)$;
\item $n\not\in \{b(v)-5,b(v)-3,b(v)-2,b(v)-1\}$;
\item if $v=9$, then $n\not\in\{5,8\}$;
\item $m\not\in\{1,2,4\}$;
\item if $n>b(v)-(v-1)/2$, then $m=v$; and
\item if $m<v$ and $n=b(m)$, then $v\geq 2m+1$.
\end{enumerate}
\end{lemma}

\noindent {\em Proof.}
To see that (i) holds, note that there exists a pair of STS$(v)$ whose intersection is $\I$, such that
$(\supp(\I),\I)$ is a partial triple system of order $m$ having $n$ blocks. Hence, $n\leq D(2,3,m)$.
That $m/3\leq n$ follows easily from the observation that the maximum number of points underlying
$n$ blocks of size three is $3n$.

Theorem \ref{LR} together with the observation that $\Int_\Phi(\{3\},1^v)|_2=\Int_{|\cdot|}(\{3\},1^v)$
gives conditions (ii) and (iii).

The observation that $\Int_\Phi(\{3\},1^v)|_1=\Int_{|\supp(\cdot)|}(\{3\},1^v)$
and that the support of the set of blocks of a partial triple system can never contain one, two, or four
points gives condition (iv).

For condition (v), observe that each point of an STS$(v)$ lies in $(v-1)/2$ blocks. So unless we have less
than $b(v)-(v-1)/2$ blocks, we cannot drive the support down to less than $v$ points.

Any partial triple system of order $m$ and $n=b(m)$ blocks is an STS$(m)$. Theorem \ref{DW}
then implies that $v\geq 2m+1$, from which condition (vi) follows.
$\Box$ \\

We call the set
\begin{eqnarray*}
{\sf A}(v) & = & \{(m,n)\in\mathbb{Z}^2_{\geq 0} : \text{$(m,n)$ satisfies conditions (i)-(vi) of Lemma
\ref{admissible}}\}
\end{eqnarray*}
the {\em admissible region} for $\Int_\Phi(\{3\},1^v)$.

\begin{lemma}
\label{up}
Let $v\equiv 1$ or $3\pmod{6}$. Then $|\Int_\Phi(\{3\},1^v)|\leq (1+o(1))\frac{1}{18} v^3$.
\end{lemma}

\noindent {\em Proof.}
\begin{eqnarray*}
|\Int_\Phi(\{3\},1^v) & \leq & |{\sf A}(v)| \\
& \leq & \sum_{m=0}^v (D(2,3,m)-m/3+1)~~~~\text{(by condition (i) of Lemma \ref{admissible})}\\
& = & \sum_{m=0}^v \left(\frac{1}{6} m^2 +o(m^2)\right) \\
& = & (1+o(1))\frac{1}{18} v^3.
\end{eqnarray*}
$\Box$ \\

Existing results on the intersection of STS only determine that a negligible portion of ${\sf A}(v)$
belongs to $\Int_\Phi(\{3\},1^v)$. In particular, Theorem \ref{LR} only implies that
$|\Int_\Phi(\{3\},1^v)|\geq (1+o(1))\frac{1}{6}v^2$. The main result of this paper is the following.\\

\noindent {\bf Main Theorem}~~{\em For $v\equiv 1$ or $3\pmod{6}$,
$|\Int_\Phi(\{3\},1^v)|=\Theta(v^3)$.}

\section{Proof of the Main Theorem}

Our main tool is Wilson's Fundamental Construction for GDDs \cite{Wilson:1972a}.

\begin{center}
\begin{tabular}{| l l|}
\hline
\multicolumn{2}{|c|}{\underline{Wilson's Fundamental Construction}}\\
Input: & (master) GDD $\D=(X,\G,\A)$; \\
           & weight function $\omega:X\rightarrow\mathbb{Z}_{\geq 0}$; \\
           & (ingredient) $K$-GDD $\D_A=(X_A,\G_A,\B_A)$ of type $[\omega(a):a\in A]$,\\
           & for each block $A\in\A$, where \\
           & ~~~~$X_A=\cup_{a\in A}\{\{a\} \times \{1,\ldots,\omega(a)\}\}$ and \\
           & ~~~~$\G_A=\{\{a\}\times\{1,\ldots,\omega(a)\}:a\in A\}$. \\
Output: & $K$-GDD $\D^*=(X^*,\G^*,\A^*)$ of type $[\sum_{x\in G}\omega(x):G\in\G]$, where \\
              & ~~~~$X^*=\cup_{x\in X}(\{x\}\times\{1,\ldots,\omega(x)\})$, \\
              & ~~~~$\G^*=\{\cup_{x\in G}(\{x\}\times\{1,\ldots,\omega(x)\}) : G\in\G\}$, and \\
              & ~~~~$\A^*=\cup_{A\in\A}\B_A$. \\
Notation: & $\D^*={\rm WFC}(\D,\omega,\{\D_A:A\in\A\})$. \\
Note: & By convention, for $x\in X$, $\{x\}\times\{1,\ldots,\omega(x)\}=\varnothing$ if $\omega(x)=0$.\\
\hline
\multicolumn{2}{c}{}\\
\end{tabular}
\end{center}

The master GDDs we use are the class of $\{4\}$-GDDs of type $1^ut^1$, existence for which has
been settled by Rees and Stinson \cite{ReesStinson:1989b}.

\begin{theorem}[Rees and Stinson \cite{ReesStinson:1989b}]
\label{RS}
There exists a $\{4\}$-GDD of type $1^ut^1$ whenever $u\geq 2t+1$ and 
\begin{enumerate}[{\rm (i)}]
\item $u\equiv 0$ or $3\pmod{12}$ and
$t\equiv 1$ or $7\pmod{12}$; or
\item $u\equiv 0$ or $9\pmod{12}$ and $t\equiv 4$ or $10\pmod{12}$.
\end{enumerate}
\end{theorem}

We call the distinguished group of size $t$ in a $\{4\}$-GDD of type $1^ut^1$, the {\em hole}.

\begin{lemma}
The number of blocks in a $\{4\}$-GDD of type $1^ut^1$ that are disjoint from the hole is
$u(u-2t-1)/12$.
\end{lemma}

\noindent {\em Proof.}
The number of blocks that have non-empty intersection with the hole is easily seen to be $ut/3$.
The total number of blocks in the GDD is $\left({u+t\choose 2}-{t\choose 2}\right)/6$. Hence, the number
of blocks disjoint from the hole is $\left({u+t\choose 2}-{t\choose 2}\right)/6-ut/3=u(u-2t-1)/12$.
$\Box$ \\

We also make use of the following result of Butler and Hoffman \cite{ButlerHoffman:1992}.

\begin{theorem}[Butler and Hoffman \cite{ButlerHoffman:1992}]
\label{BH}
Let $g$ and $t$ be positive integers such that $t\geq 3$, $g^2{t\choose 2}\equiv 0\pmod{3}$, and
$g(t-1)\equiv 0\pmod{2}$. Let $b(g^t)=g^2t(t-1)/6$ and denote by $I(g^t)=\{0,\ldots,b(g^t)\}\setminus
\{b(g^t)-5,b(g^t)-3,b(g^t)-2,b(g^t)-1\}$. Then $\Int_{|\cdot|}(\{3\},g^t)=I(g^t)$, except that
\begin{enumerate}[{\rm (i)}]
\item $\Int_{|\cdot|}(\{3\},1^9)=I(1^9)\setminus\{5,8\}$;
\item $\Int_{|\cdot|}(\{3\},2^4)=I(2^4)\setminus\{1,4\}$;
\item $\Int_{|\cdot|}(\{3\},3^3)=I(3^3)\setminus\{1,2,5\}$; and
\item $\Int_{|\cdot|}(\{3\},4^3)=I(4^3)\setminus\{5,7,10\}$.
\end{enumerate}
\end{theorem}

Let $\D=(X,\G,\A)$ be a $\{4\}$-GDD of type $1^ut^1$, with $\G=\{G_1,\ldots,G_{u+1}\}$, where
\begin{eqnarray*}
G_i & = & \begin{cases}
\{x_i\},&\text{if $1\leq i\leq u$; and} \\
\{x_{u+1},\ldots,x_{u+t}\},&\text{if $i=u+1$.}
\end{cases}
\end{eqnarray*}
For $\alpha,\beta\geq 0$ such that $\alpha+\beta\leq t$, define the following weight function:
\begin{eqnarray*}
\omega_{\alpha,\beta}(x) & = & \begin{cases}
2,&\text{if $x\in\{x_1,\ldots,x_{u+\alpha}\}$;} \\
4,&\text{if $x\in\{x_{u+\alpha+1},\ldots,x_{u+\alpha+\beta}\}$; and} \\
0,&\text{if $x\in\{x_{u+\alpha+\beta+1},\ldots,x_{u+t}\}$.}
\end{cases}
\end{eqnarray*}
We use Wilson's Fundamental Construction with $\D$ as master GDD and $\omega_{\alpha,\beta}$
as weight function. The required ingredient GDDs are $\{3\}$-GDDs of type $2^4$, type
$2^34^1$, and type $2^3$. The following are results on the intersections of these ingredient GDDs.

\begin{lemma}
\label{small}
The following hold:
\begin{enumerate}[{\rm (i)}]
\item $\Int_{|\cdot|}(\{3\},2^3)=\{0,4\}$;
\item $\Int_{|\cdot|}(\{3\},2^4)=\{0,2,8\}$; and
\item $0\in \Int_{|\cdot|}(\{3\},2^3 4^1)$;
\end{enumerate}
\end{lemma}

\noindent {\em Proof.}
(i) and (ii) follow from Theorem \ref{BH}. The existence of a pair of disjoint $\{3\}$-GDDs of
type $2^3 4^1$ is given in \cite{SchellenbergStinson:1989}, proving (iii).
$\Box$ \\

For a set $A$ of integers, we denote by $A+A$ the set $\{a+b:a,b\in A\}$, and denote by
$\sum_{i=1}^s A$ the set $A+\cdots+A$ ($s$-fold sum).

\begin{lemma}
\label{interval}
Let $s\geq 2$. Then $|\Int_{|\cdot|}(\{3\},2^4)|=4s-2$.
\end{lemma}

\noindent {\em Proof.}
We have $\sum_{i=1}^s \Int_{|\cdot|}(\{3\},2^4)=\{0,2,4,\ldots,8s\}\setminus\{8s-10,8s-4,8s-2\}$.
$\Box$

\begin{lemma}
\label{chooseb}
Let $m\geq 0$ and $r\in\{1,3,7,9,13,15,19,21\}$. Then for every $a\geq\left\lceil (m+3)/2\right\rceil$,
there exists an integer $b\geq 0$ such that the following inequalities are all satisfied:
\begin{enumerate}[{\rm (i)}]
\item $24(m-a)+r\leq 4(12b+1)+1$;
\item $12a\geq 2(12b+1)+1$.
\end{enumerate}
\end{lemma}

\noindent {\em Proof.}
Let $b$ be the smallest integer such that inequality (i) holds. Then
\begin{eqnarray*}
b=\left\lceil \frac{24(m-a)+r-5}{48}\right\rceil \geq \left\lceil \frac{r-5}{48}\right\rceil \geq
\left\lceil -\frac{4}{48}\right\rceil = 0.
\end{eqnarray*}
It remains to show that inequality (ii) holds. To see that this is the case, observe that
\begin{eqnarray*}
2(12b+1)+1 & = & 24b+3 \\
& = & 24\left\lceil \frac{24(m-1)+r-5}{48}\right\rceil + 3\\
& = & 24 \left\lceil \frac{m-a}{2} + \frac{r-5}{48}\right\rceil + 3\\
& \leq & 24\left\lceil \frac{m-a}{2} + \frac{16}{48}\right\rceil+ 3\\
& \leq & 24\left(\frac{m-a}{2}+1\right)+3 \\
& = & 12(m-a)+27 \\
& \leq & 12(2a-3-a)+27 \\
& = & 12a-9\\
& \leq & 12a.
\end{eqnarray*}
$\Box$ \\

Let $m\geq 0$ and $r\in\{1,3,7,9,13,15,19,21\}$. We now construct a pair of STS$(24m+r)$
via Wilson's Fundamental Construction. Write $24m+r$ as $24a+24(m-1)+r$, where
$a\geq\left\lceil (m+3)/2\right\rceil$. Choose $b$ to be the smallest non-negative integer so that
the inequalities $24(m-a)+r\leq 4(12b+1)+1$ and $12a\geq 2(12b+1)+1$ are both satisfied. Such
a $b$ exists by Lemma \ref{chooseb}. Choose also non-negative integers $\alpha$ and $\beta$
such that $0\leq\alpha+\beta\leq 12b+1$ and $2\alpha+4\beta+1=24(m-a)+r$. This is always possible
because $\{2\alpha+4\beta+1:0\leq \alpha+\beta\leq 12b+1\}=\{1,3,5,\ldots,4(12b+1)+1\}\ni 24(m-a)+r$.

Now, take $\D=\{X,\G,\A)$ to be a $\{4\}$-GDD of type $1^{12a}(12b+1)^1$. The existence of such
a GDD is implied by the inequality $12a\geq 2(12b+1)+1$ (via Theorem \ref{RS}). Let
${\sf G}={\rm WFC}(\D,\omega_{\alpha,\beta},\{\D_A:A\in\A\})$ and
${\sf G}'={\rm WFC}(\D,\omega_{\alpha,\beta},\{\D'_A:A\in\A\})$, where
$\D_A$ and $\D'_A$ are a pair of
\begin{enumerate}[{\rm (i)}]
\item $\{3\}$-GDDs of type $2^4$ intersecting in $\mu_A$ blocks, if $A$ is disjoint from the hole of
$\D$;
\item disjoint $\{3\}$-GDDs of type $2^4$, if $A$ contains a point of weight two from the hole of $\D$;
\item disjoint $\{3\}$-GDDs of type $2^34^1$, if $A$ contains a point of weight four from the hole of $\D$;
and
\item disjoint $\{3\}$-GDDs of type $2^3$, if $A$ contains a point of weight zero from the hole of $\D$.
\end{enumerate}
Such ingredient GDDs all exist by Lemma \ref{small}.

It is clear from the description of Wilson's Fundamental Construction that $\sf G$ and ${\sf G}'$
are two $\{3\}$-GDDs of type $2^{12a}(2\alpha+4\beta)^1$, where $0\leq\alpha+\beta\leq 12b+1$,
intersecting in $\sum_{A\in\A}\mu_A$ blocks, whose support is disjoint from the group of size
$2\alpha+4\beta$. Now add a point to each of $\sf G$ and ${\sf G}'$ to obtain
$\{3,2\alpha+4\beta+1\}$-GDDs of type $1^{24a+2\alpha+4\beta+1}$ with exactly one block of
size $2\alpha+4\beta+1$. Replace the block of size $2\alpha+4\beta+1=24(m-a)+r$ in each of these
GDDs with the respective blocks from a pair of disjoint STS$(24(m-a)+r)$ (which exists by Theorem
\ref{LR}). The result is a pair of STS$(24m+r)$ intersecting in $12a+\sum_{A\in\A}\mu_A$ blocks
whose support contains exactly $24a+1$ points.

By varying $a$ and $\mu_A$, this shows that
\begin{eqnarray}
\label{ineq}
|\Int_\Phi(\{3\},1^{24m+r})| & \geq & \sum_{a=\left\lceil (m+3)/2\right\rceil}^m
\left| \sum_{i=1}^{a(12a-2(12b+1)-1)} \Int_{|\cdot|}(\{3\},2^4)\right|.
\end{eqnarray}

\begin{lemma}
\label{a}
For $a\geq\left\lceil (m+3)/2\right\rceil$, we have $12a-2(12b+1)-1\geq 24a-12m-27$.
\end{lemma}

\noindent {\em Proof.}
By the proof of Lemma \ref{chooseb}, we know that
$b\leq\left\lceil \frac{24(m-a)+r-5}{48}\right\rceil\leq 1+(m-a)/2$. So,
\begin{eqnarray*}
12a-2(12b+1)-1 & = & 12a-24b-3 \\
& \geq & 12a-24(1+(m-a)/2)-3 \\
& = & 24a-12m-27.
\end{eqnarray*}
$\Box$ \\

By Lemma \ref{a}, inequality (\ref{ineq}) implies
\begin{eqnarray*}
|\Int_\Phi(\{3\},1^{24m+r})| & \geq & \sum_{a=\left\lceil (m+3)/2\right\rceil}^m
\left| \sum_{i=1}^{a(24a-12m-27)} \Int_{|\cdot|}(\{3\},2^4)\right| \\
& = & \sum_{a=\left\lceil (m+3)/2\right\rceil}^m (4a(24a-12m-27)-2)~~~~
\text{(via Lemma \ref{interval})} \\
& = & \left(96\sum_{a=\left\lceil (m+3)/2\right\rceil}^m a^2\right) -
\left(4(12m+27)\sum_{a=\left\lceil (m+3)/2\right\rceil}^m a\right) - \\
& & \left(\sum_{a=\left\lceil (m+3)/2\right\rceil}^m 2\right) \\
& \geq & (28m^3 +o(m^3))-(18m^3+o(m^3))-o(m^3) \\
& = & 10m^3 +o(m^3).
\end{eqnarray*}
This shows that for $v\equiv 1$ or $3\pmod{6}$,
\begin{eqnarray*}
|\Int_\Phi(\{3\},1^v)| & \geq & (1+o(1))\frac{5}{6912} v^3,
\end{eqnarray*}
which together with Lemma \ref{up} proves the Main Theorem.

\section{Conclusion}

In this paper, we initiated the study on the fine intersection problem for STS.
We established that $|\Int_\Phi(\{3\},1^v)|=\Theta(v^3)$ for $v\equiv 1$ or $3\pmod{6}$. There
remain many interesting unsolved problems:
\begin{enumerate}[{\rm (i)}]
\item What is the exact asymptotics of $|\Int_\Phi(\{3\},1^v)|$? There remains a wide gap
between our lower and upper bounds on $|\Int_\Phi(\{3\},1^v)|$. We think the upper bound is
probably the truth and make the conjecture that $|\Int_\Phi(\{3\},1^v)|=(1+o(1))\frac{1}{18}v^3$.
\item Determine completely the set $\Int_\Phi(\{3\},1^v)$. This problem is probably very difficult.
\item What is the number of non-isomorphic partial triple systems that can underly the
intersection of two STS$(v)$?
\item Determine all non-isomorphic partial triple systems that can underly the intersection of
two STS$(v)$.
\end{enumerate}

Intersection problems for STS remain well alive three decades after the seminal paper of
Lindner and Rosa \cite{LindnerRosa:1975a},

\end{document}